\documentclass[12 pt]{article}
\usepackage{amsmath}
\usepackage{amsthm}
\usepackage{graphicx}
\usepackage{amsmath,bbm,amssymb}
\usepackage{bbold}

\newcommand{\ceil}[1]{\lceil {#1} \rceil}

\newtheorem{theorem}{Theorem}[section]
\newtheorem{prop}[theorem]{Proposition}
\newtheorem{lemma}[theorem]{Lemma}

\theoremstyle{remark}
\newtheorem*{remark}{Remark}

\theoremstyle{definition}

\title{Spherical Asymptotics for the \\ Rotor-Router Model in $\Z^d$}
\author{Lionel Levine \\ Yuval Peres}
\date{March 13th, 2005}

\DeclareSymbolFont{AMSb}{U}{msb}{m}{n}
\DeclareMathSymbol{\C}{\mathbin}{AMSb}{"43}
\DeclareMathSymbol{\EE}{\mathbin}{AMSb}{"45}
\DeclareMathSymbol{\N}{\mathbin}{AMSb}{"4E}
\DeclareMathSymbol{\PP}{\mathbin}{AMSb}{"50}
\DeclareMathSymbol{\Q}{\mathbin}{AMSb}{"51}
\DeclareMathSymbol{\R}{\mathbin}{AMSb}{"52}
\DeclareMathSymbol{\Z}{\mathbin}{AMSb}{"5A}

\begin{document}
\maketitle

\begin{abstract} 
The rotor-router model is a deterministic analogue of random walk invented by Jim Propp.  It can be used to define a
deterministic aggregation model analogous to internal diffusion limited aggregation.  We prove an isoperimetric inequality for the exit
time of simple random walk from a finite region in $\Z^d$, and use this to prove that the shape of the rotor-router aggregation model
in $\Z^d$, suitably rescaled, converges to a Euclidean ball in $\R^d$.  
\end{abstract}

\section{Introduction}
Given a finite region $A \subset \Z^d$, let $A'$ be the (random) region obtained by starting a random walk at the origin, stopping the walk when it first exits $A$, and 
adjoining the endpoint of the walk to $A$. Internal diffusion limited aggregation (``internal DLA'') is the growth model obtained by iterating this procedure starting from 
the set containing only the origin: $A_1 = \{o\}$, $A_n = (A_{n-1})'$.  Lawler et al.\ \cite{LBG} showed that the region $A_n$, rescaled by a factor of $n^{1/d}$, 
converges with probability one to a Euclidean ball in $\R^d$ as $n \rightarrow \infty$.  Lawler \cite{Lawler} estimated the rate of convergence.

Jim Propp has proposed the following deterministic analogue of internal DLA in two dimensions.  At each site $x \in A$ is a ``rotor'' pointing North, East, South or West.  
A particle is placed at the origin and performs {\it rotor-router walk} until it exits the region $A$: during each time step, the rotor at the particle's current location 
is rotated clockwise by $90$ degrees, and the particle takes a step in the direction of the newly rotated rotor.  The intent of this rule is to simulate the first-order
properties of random walk by forcing each site to route approximately equal numbers of particles to each of the four neighboring sites.  When the particle reaches a point
not in $A$, that point is adjoined to the region and the procedure is iterated to obtain a sequence of regions $A_n$.  For example, if 
all rotors are initially pointing north, the sequence will be begin $A_1 = \{o\}$, $A_2 = \{o,(1,0)\}$, $A_3 = \{o,(1,0),(0,-1)\}$, etc.

\begin{figure}
\label{blobs}
\centering
\includegraphics[scale=.25]{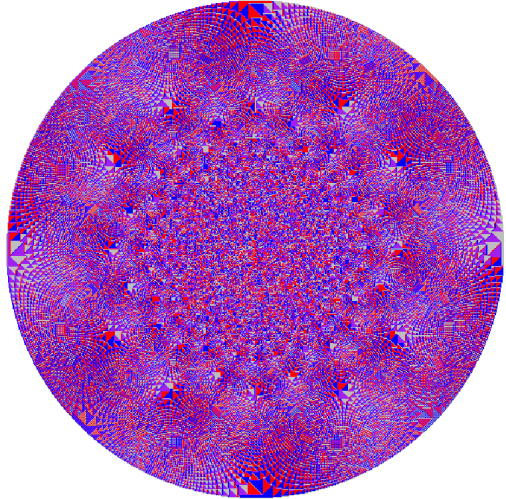}
\caption{Rotor-router aggregate of $270,000$ particles.}
\end{figure}

There has been considerable recent interest in obtaining a shape theorem for the rotor-router model analogous to that for internal DLA \cite{Kleber, thesis}.  Much of this 
interest has been driven by simulations in two dimensions, which indicate that the regions $A_n$ are extraordinarily close to circular (Figure~\ref{blobs}).  
Despite the impressive evidence for circularity, very little progress has been made until now in the way of rigorous results.  In one dimension, with rotors alternately 
pointing left and right, the dynamics of the model are simple enough to analyze explicitly; in this case the first author has shown \cite{thesis} that the deviation from a 
ball (symmetric interval) is bounded independent of $n$.  In addition, various modifications and extensions of the one-dimensional model are amenable to explicit analysis, 
and analogous shape theorems are known in some of these cases \cite{Kleber, thesis}.  In two dimensions, the first author has shown \cite{thesis} that the region $A_n$ 
contains a disc of radius proportional to $n^{1/4}$.  In higher dimensions, the model can be defined analogously by repeatedly cycling the rotors through an ordering of the 
$2d$ cardinal directions in $\Z^d$; until now nothing was known about the shape for $d \geq 3$.

Denote by $R ~\Delta~ S$ the symmetric difference of sets $R$ and $S$.  For a region $A \subset \Z^d$, we write $A^{\Box}$ for the union of closed unit cubes in $\R^d$ 
centered at the points of $A$.  We write $\mathcal{L}$ for $d$-dimensional Lebesgue measure.  As a special case of our main result, Theorem~\ref{main}, we obtain the 
following.

\begin{theorem}
\label{intro}
Let $(A_n)_{n\geq 1}$ be rotor-router aggregation in $\Z^d$, starting from any initial configuration of rotors.  Then as $n \rightarrow \infty$
	\[ \mathcal{L}( n^{-1/d} A_n^\Box ~\Delta~ B ) \rightarrow 0, \]
where $B$ is the ball of unit volume centered at the origin in $\R^d$.
\end{theorem}

In Theorem~\ref{main} we prove this result in the more general setting of arbitrary rotor stacks of bounded ``discrepancy;'' see section~2 for details.  We also give an 
explicit bound on the rate of convergence.  Here we should emphasize that much work remains to be done if one hopes to explain the 
almost perfect circularity found in Figure~\ref{blobs}.  The form of convergence in Theorem~\ref{intro} is not as strong as the 
convergence in the shape theorems for internal DLA \cite{LBG, Lawler}.  In particular, Theorem~\ref{intro} does not preclude the 
formation of long tendrils, or of ``holes'' close to the origin, provided that the volume of these features is negligible compared to 
$n$.  So we hope that this represents only the beginning of attempts to understand the shape of the model.
 
A major component of the proof is an isoperimetric inequality for the expected exit time of random walk from a region in $\Z^d$.  Because of its intrinsic interest and 
possible utility in other applications, we state it here.  Given $x \in A \subset \Z^d$, let $e_x(A)$ be the expected time for simple random walk started at $x$ to first 
leave the region $A$.  We denote by $\omega_d$ the volume of the ball of unit radius in $\R^d$.  The following result shows that $e_x(A)$ is asymptotically maximized among 
all regions of a given size when $A$ is a ball.

\begin{theorem}
If $x \in A \subset \Z^d$ and $|A|=n$, then for sufficiently small $\epsilon>0$
	\[ e_x(A) \leq e_o(B_n)(1+O(n^{-\epsilon})), \]
where $B_n$ is the lattice ball of radius $(n/\omega_d)^{1/d}$:
	\begin{equation} \label{B_n} B_n = \{ y \in \Z^d ~|~ \omega_d ||y||^d < n \}. \end{equation}
\end{theorem}

Here $o$ denotes the origin in $\Z^d$ and $||\cdot||$ the Euclidean norm.  In Theorem~\ref{isoperimetric} we give an explicit bound for the exponent in the error term.

The proof of Theorem~\ref{intro} also relies on a special ``abelian'' property of internal DLA.  Beginning with a configuration of 
finitely many particles in $\Z^d$, suppose that at each time step we choose a site occupied by more than one particle and move one 
particle at that site by one step according to simple random walk.  After a succession of such choices, the process ends when no site 
is occupied by more than one particle.  Diaconis and Fulton \cite{DF} discovered the remarkable fact that both the distribution of the 
shape of the resulting cluster of occupied sites and the distribution of the number of time steps taken to reach the final state are 
independent of the choices made along the way.  Lawler et al.\ \cite{LBG} (sec.\ 6) found a particularly simple proof by labeling the 
particles and moving only the highest-labeled particle at a given site.

\section{Convergence to a Ball}

For $x$, $y \in \Z^d$ we write $x \sim y$ if $||x-y||=1$.  By a ``region'' $A \subset \Z^d$ we will always mean a finite region.  We 
write $|A|$ for the cardinality of $A$.  The boundary of $A$ is the region
	\[ \partial A = \{ x \in A^c ~|~ x \sim y \text{ for some } y \in A \}. \]
Simple random walk in $\Z^d$ will be denoted by $X_0, X_1, \ldots$; the probability and expectation operators for walk started at $X_0=x$ will be denoted $\PP_x$ and 
$\EE_x$.  Given a set $S \subset \Z^d$ we write
	\[ T_S = \inf \{ j \geq 0 ~|~ X_j \in S \} \]
for the first hitting time of $S$.  If $A$ is a region in $\Z^d$ and $x\in A$, we denote by $e_x(A)$ the expected time for simple random walk started at $x$ to exit the set 
$A$:
	\[ e_x(A) = \EE_x T_{\partial A}. \]
Note that $e_x(A) = 0$ for $x \notin A$.  We write
	\begin{equation} \label{bar} \bar{e}(A) = \sup_x e_x(A). \end{equation}

We study the following mild generalization of the rotor-router model in $\Z^d$.  Fix a positive constant $D$, the {\it discrepancy}.  At each site $x \in \Z^d$ is an 
infinite stack of rotors $r_1, r_2, \ldots$ each pointing in one of the $2d$ cardinal directions.  On the $i$-th visit to the site $x$, the particle travels in 
direction $r_i$.  We require that for any direction $\varepsilon$ and any positive integer $m$,
	\begin{equation} \label{discrepancy} \left| \# \{ i\leq m ~|~ r_i = \varepsilon \} - \frac{m}{2d} \right| \leq D. \end{equation}
Observe that the original rotor-router model with cyclically repeating rotors satisfies this condition with discrepancy $D=1$.

We will phrase our result in terms of the quantity
	\[ \phi(n) = \sup_{A \subset \Z^d, ~|A|=n} (\bar{e}(A) - e_o(B_n)). \]
Let 
	\[ \Phi(n) = \sum_{j=1}^n \phi(j). \]
In Theorem~\ref{isoperimetric} we show that $\phi(n) = O(n^{2/d-\epsilon})$ for sufficiently small $\epsilon$.  

Write $\mathcal{L}$ for Lebesgue measure in $\R^d$.  For $A \subset \Z^d$ we write
	\[ A^\Box = A + \left[ -\frac12, \frac12 \right]^d \subset \R^d. \]
For $R \subset \R^d$ and $\lambda \in \R$, we write
	\[ \lambda R = \{ \lambda x ~|~ x\in R \}.\]
Our main result is the following.

\begin{theorem}
\label{main}
Let $(A_n)_{n \geq 1}$ be rotor-router aggregation in $\Z^d$ using any configuration of rotor stacks with discrepancy at most $D$.  Then
	\[ \mathcal{L}(n^{-1/d}A_n^\Box ~\Delta~ B) = C n^{-1/2-1/d} \Phi(n)^{1/2} + O(D^{1-1/8d} n^{-1/2d}) \rightarrow 0 \]
as $n \rightarrow \infty$, where $B$ is the ball of unit volume centered at the origin in $\R^d$, and $C=C(d)$ is a constant 
independent of $n$ and $D$.
\end{theorem}

\begin{remark} 
In two dimensions, Theorem~\ref{isoperimetric} implies that $\Phi(n) = O(n^{5/3})$, hence
	\[ \mathcal{L}(n^{-1/d} A_n^\Box ~\Delta~ B) = O(n^{-1/6}) + O(D^{15/16} n^{-1/4}). \]
For $d\geq 3$, Theorem~\ref{isoperimetric} gives
	\begin{equation} \label{Phi} \Phi(n) = O(n^{1+2/d-\frac{2^{-d}}{2 d^2 \log 3}}), \end{equation}
hence
	\[ \mathcal{L}(n^{-1/d} A_n^\Box ~\Delta~ B) = O(n^{-\frac{2^{-d}}{4 d^2 \log 3}}) + O(D^{1-1/8d} n^{-1/2d}). \]
\end{remark}

Given a finite region $A \subset \Z^d$, define
	\[ \psi(A) = \sum_{x \in A} ||x||^2. \]
Among regions $A \subset \Z^d$ with $|A|=|B_n|$, the quantity $\psi(A)$ is minimized when $A=B_n$.  The idea of the proof of Theorem~\ref{main} is to show that $\psi(A_n)$ 
cannot be much larger than $\psi(B_n)$ (Proposition~\ref{mainprop}).  To estimate $\psi(B_n)$, notice that
	\begin{equation} \label{balls} B((n/\omega_d)^{1/d}-d) \subset B_n^\Box \subset B((n/\omega_d)^{1/d}+d), \end{equation}
where $B(r)$ is the ball of radius $r$ centered at the origin in $\R^d$.  From this we obtain
	\begin{equation} \label{psi} \psi(B_n) = \frac{d}{d+2} \omega_d^{-2/d} n^{1+2/d} + O(n^{1+1/d}). \end{equation}

\begin{prop}
\label{mainprop}
$\psi(A_n) = \psi(B_n) + \Phi(n) + O(D^{2-1/4d} n^{1+1/d})$.
\end{prop}

To prove Proposition~\ref{mainprop}, we first relate the quantity $\psi(A_n)$ to the total number of steps $T_n$ taken by the 
rotor-router walks of the first $n$ particles. We will make use of the identity
	\begin{equation} \label{r^2} \Delta_x ||x||^2 = \frac{1}{2d} \sum_{i=1}^d ((x_i-1)^2 - 2x_i^2 + (x_i+1)^2) = 1. \end{equation}
where $\Delta_x f$ denotes the discrete Laplacian of the 
function $f$ at the point $x$:
	\[ \Delta_x f = \frac{1}{2d} \sum_{y \sim x} f(y) - f(x) = \EE_x f(X_1) - f(x). \]

\begin{lemma}
\label{time}
$\psi(A_n) \leq T_n + 8 \sqrt{d} D \sum_{x \in A_n} ||x|| + 4dDn$.
\end{lemma}

\begin{proof}
Given a finite set of particles at locations $x(1), \ldots, x(n) \in \Z^d$, define the quadratic weight of the configuration to be
	\[ w = w(x(1), \ldots, x(n)) = \sum_{i=1}^n ||x(i)||^2. \]
At any given time during rotor-router aggregation, each site $x$ has routed an equal number of particles to each of its neighbors, plus an error of at most $2D$ extra 
routings to each neighbor $y \sim x$. By (\ref{r^2}) it follows that the net effect of the first $m$ routings from a site $x$ is to 
increase the total weight by $m$ plus an error of at most $2D \sum_{y \sim x} \left| ||x||^2 - ||y||^2 \right|$.

Starting with $n$ particles at the origin, let the particles perform rotor-router aggregation one at a time until all of $A_n$ is occupied.  This process involves a total 
of $T_n$ routings, so the net change in weight is $T_n$ plus an error of at most
	\begin{eqnarray*} 2D \sum_{x \sim y \in A_n} \left| ||x||^2 - ||y||^2 \right|
				&=& 2D \sum_{x \in A_n} \sum_{i=1}^d ( |2x_i + 1| + |2x_i - 1| )  \\
				&\leq& 2D \sum_{x \in A_n} \sum_{i=1}^d (4|x_i|+2). \\
				&\leq& 8 \sqrt{d} D \sum_{x \in A_n} ||x|| + 4dDn, \end{eqnarray*}
where the last inequality is Cauchy-Schwarz.  The result now follows from the fact that the initial configuration has weight zero and the final configuration has weight 
$\psi(A_n)$.
\end{proof}

We next relate the quantity $T_n$ to the expected exit time $e_x(A_n)$.  For any region $A$ we have the Laplacian identity
	\begin{equation} \label{laplacian} \Delta_x e_x(A) =  -1, \qquad x \in A, \end{equation}
We will need an estimate for the exit time $e_o(B_n)$ from the discrete ball $B_n$.  By stopping the bounded martingale $||X_t||^2 - t$ 
at the time $T=T_{\partial B_n}$ when the walk exits $B_n$, we obtain
	\begin{equation} \label{martingale} e_o(B_n) = \EE_{o}T = \EE_{o}||X_T||^2 = (n/\omega_d)^{2/d} + O(n^{1/d}). \end{equation}
This asymptotic result for random walk becomes exact for Brownian motion when $B_n$ is replaced by a ball of radius 
$(n/\omega_d)^{1/d}$ in $\R^d$.

\begin{lemma} \label{cauchyschwarz} If $|A|=n$, then $\sum_{x \sim y \in A} |e_x(A) - e_y(A)| = O(n^{1+1/d})$. \end{lemma}
\begin{proof}
By Cauchy-Schwarz,
	\[ \left( \sum_{x \sim y \in A} |e_x(A) - e_y(A)| \right)^2 \leq 4n \sum_{x \sim y \in A} (e_x(A)-e_y(A))^2. \]
To bound the latter sum, the fact that $e_x(A)=0$ for $x \in \partial A$, together with (\ref{laplacian}) and Theorem~\ref{isoperimetric}, implies
\begin{eqnarray*}  \sum_{x \sim y \in A \cup \partial A} (e_x(A) - e_y(A))^2 &=&  2 \sum_{x \sim y \in A\cup\partial A} e_x(A) (e_x(A) - e_y(A)) \\
			 &=& 2\sum_{x \in A} e_x(A) \sum_{y \sim x} (e_x(A) - e_y(A)) \\
			 &=& 2\sum_{x \in A} e_x(A) (-2d\Delta_x e_x(A)) \\
			 &=& 4d\sum_{x \in A} e_x(A) \\
			 &=& O(n^{1+2/d}). \qed \end{eqnarray*}
\renewcommand{\qedsymbol}{}
\end{proof}

\begin{lemma} \label{exit} $T_n = n e_o(A_n) - \sum_{x \in A_n} e_x(A_n) + O(D n^{1+1/d})$. \end{lemma}

\begin{proof}
Given a finite set of particles at locations $x(1), \ldots, x(n) \in A_n$, define the exit weight of the configuration to be 
	\[ \eta = \eta(x(1), \ldots, x(n)) = \sum_{j=1}^n e_{x(j)}(A_n). \]
By (\ref{discrepancy}) and (\ref{laplacian}), the net effect of the first $m$ routings from a site $x$ is to decrease the total 
weight $\eta$ by $m$, plus an error of at most $2D \sum_{y \sim x} |e_x(A_n) - e_y(A_n)|$.  Beginning with $n$ particles at the origin 
and ending when $A_n$ is completely occupied, the total decrease in weight thus differs from $T_n$ by at most $2D \sum_{x \sim y \in 
A_n} |e_x(A_n)-e_y(A_n)|$.  The result now follows from Lemma~\ref{cauchyschwarz}.
\end{proof}

\begin{lemma} \label{T_n} $T_n \leq \frac{d}{d+2} \omega_d^{-2/d} n^{1+2/d} + \Phi(n) + O(D n^{1+1/d})$. \end{lemma}
 
\begin{proof} It is here that we will use the abelian property of internal DLA mentioned in the introduction.  Beginning with $n$ 
particles at the origin, let each particle $p_j$ in turn perform simple random walk until it either exits the region $A_n$ or visits a 
site not visited by any of $p_1, \ldots, p_{j-1}$.  The expected time taken by this procedure is at least $n e_o(A_n) - \sum_{x\in 
A_n}e_x(A_n)$.  Letting the particles which were stopped upon exiting $A_n$ continue walking until they also reach unoccupied sites, 
the abelian property implies that
	\[ n e_o(A_n) - \sum_{x\in A_n}e_x(A_n) \leq T_n^{\text{IDLA}} \]  
where $T_n^{\text{IDLA}}$ is the expected number of steps taken by the first $n$ particles in internal DLA.  By (\ref{martingale}) we have
	\begin{eqnarray*} T_n^{\text{IDLA}} &\leq& \sum_{j=1}^n (e_o(B_j) + \phi(j)) \\
				&\leq& \sum_{j=1}^n ((j/\omega_d)^{2/d} + O(j^{1/d})) + \Phi(n) \\
				&=& \frac{d}{d+2} \omega_d^{-2/d} n^{1+2/d} + \Phi(n) + O(n^{1+1/d}). \end{eqnarray*}
The result now follows from Lemma~\ref{exit}. 
\end{proof}

We can now prove Proposition~\ref{mainprop} by means of a bootstrapping argument.

\begin{lemma} 
\label{bootstrapping} 
If $\psi(A_n) = O(D^\alpha n^\beta)$ for some $\alpha \geq 1$, $\beta \geq 1+\frac2d$, then
	\[ \psi(A_n) = O(D^{1+\alpha/2} n^{(1+\beta)/2}) + O(n^{1+2/d}). \] 
\end{lemma}

\begin{proof}
By Cauchy-Schwarz
	\begin{equation} \label{iterateme} \left( \sum_{x \in A_n} ||x|| \right)^2 \leq n \sum_{x \in A_n} ||x||^2 = n\psi(A_n) = O(D^\alpha n^{1+\beta}). \end{equation}
Lemma~\ref{time} now shows that
	\[ \psi(A_n) \leq T_n + O(D^{1+\alpha/2} n^{(1+\beta)/2}), \]
and the result follows from Lemma~\ref{T_n} and (\ref{Phi}).
\end{proof}

\begin{proof}[Proof of Proposition~\ref{mainprop}]
Since $A_n$ is connected, $||x|| \leq n$ for all $x \in A_n$, hence $\psi(A_n) = O(n^3)$.  The sequences defined by 
	\[ \alpha_0=1, \quad \alpha_{m+1}=1+\frac{\alpha_m}{2} \]
	\[ \beta_0=3, \quad \beta_{m+1}=\frac{1+\beta_m}{2} \]
have the explicit forms
	\[ \alpha_m = 2-2^{-m}, \qquad \beta_m = 1+2^{1-m}; \]
hence
	\[ \alpha_{\ceil{\log d/\log 2}} \leq 2 - \frac{1}{2d}, \quad \beta_{\ceil{\log d/\log 2}} \leq 1 + \frac2d, \]
where $\ceil{x}$ denotes the least integer $\geq x$.  By iteratively applying Lemma~\ref{bootstrapping} we obtain after $\ceil{\log d/\log 2}$ iterations
	\[ \psi(A_n) = O(D^{2-1/2d} n^{1+2/d}). \]  
Equation (\ref{iterateme}) now gives
	\[ \sum_{x \in A_n} ||x|| = O(D^{1-1/4d} n^{1+1/d}) \]
hence by Lemmas~\ref{time} and~\ref{T_n} 
	\begin{eqnarray*} \psi(A_n) &\leq& T_n + O(D^{2-1/4d} n^{1+1/d}) \\ 
				      &=& \frac{d}{d+2} \omega_d^{-d/2} n^{1+2/d} + \Phi(n) + O(D^{2-1/4d} n^{1+1/d}). \end{eqnarray*}
The result now follows from (\ref{psi}).
\end{proof}

For the proof of Theorem~\ref{main} it will be useful to rephrase equation (\ref{psi}) in terms of the radius of the ball:
	\begin{equation} \label{psirad} \psi(B_{\omega_d r^d}) = \frac{d\omega_d}{d+2} r^{d+2} + O(r^{d+1}). \end{equation}
We will also need a simple estimate for the cardinality of $B_n$.  Recall that
	\begin{equation} \label{ballsagain} B((n/\omega_d)^{1/d}-d) \subset B_n^\Box \subset B((n/\omega_d)^{1/d}+d), \end{equation}
where $B(r)$ is the ball of radius $r$ centered at the origin in $\R^d$.  It follows that
.	\begin{equation} \label{||} |B_n| = n + O(n^{1-1/d}). \end{equation}
Expressed in terms of the radius, this becomes
	\begin{equation} \label{||rad} |B_{\omega_d r^d}| = \omega_d r^d + O(r^{d-1}). \end{equation}
Better estimates exist (see e.g.\ \cite{IM}), but we will not need them.  

\begin{proof}[Proof of Theorem~\ref{main}]
We will show that
	\begin{equation} \label{enough} |A_n~\Delta~B_n| \leq C n^{1/2 - 1/d} \Phi(n)^{1/2} + O(D^{1-1/8d} n^{1-1/2d}). \end{equation}
By (\ref{ballsagain}) this implies
	\begin{eqnarray*} \mathcal{L}(n^{-1/d} A_n^\Box ~\Delta~B) &\leq& 
			\mathcal{L}(n^{-1/d} A_n^\Box ~\Delta~ n^{-1/d} B_n^\Box) + \mathcal{L}(n^{-1/d} B_n^\Box ~\Delta~ B) \\
		&\leq& \frac1n |A_n~\Delta~B_n| + O(n^{-1/d}) \\
		&\leq& C n^{-1/2-1/d} \Phi(n)^{1/2} + O(D^{1-1/8d} n^{-1/2d}), \end{eqnarray*}
which gives the theorem.

To prove (\ref{enough}), we first observe that if $|A_n~\Delta~B_n|=V$, then by (\ref{||}) we have $\psi(A_n) \geq \psi(A)$, where $A=(B_n \setminus S_-)\cup S_+$ is the 
region formed by deleting from $B_n$ an outer spherical shell $S_-$ of cardinality $V/2+O(n^{1-1/d})$ and adjoining an adjacent 
spherical shell $S_+$ of cardinality $V/2-O(n^{1-1/d})$.  The outer radius of $S_-$ and inner radius of $S_+$ are both equal to $r = 
(n/\omega_d)^{1/d}$.  Solving for the inner radius $r_-$ of $S_-$ and the outer radius $r_+$ of $S_+$ in terms of $V$, we obtain from 
(\ref{||rad})
	\[ r_{\pm} = \left(\frac{n \pm V/2+O(n^{1-1/d})}{\omega_d}\right)^{1/d} + O(1). \]
Since $t^{1+2/d}$ is a convex function of $t$, it follows that
	\begin{equation} \label{ready} r_{\pm}^{d+2} = C_0 (n \pm V/2)^{1+2/d} + O(n^{1+1/d}). \end{equation}
Equation (\ref{psirad}) yields
	\begin{eqnarray*} \psi(A_n) - \psi(B_n) &\geq& \psi \left( (B_{\omega_d r_+^d} \setminus B_{\omega_d r^d} ) \cup 
									B_{\omega_d r_-^d} \right) - \psi(B_n) \\
						&=& C_1 [ r_+^{d+2} - 2r^{d+2} + r_-^{d+2}] + O(n^{1+1/d}). \end{eqnarray*}
Applying (\ref{ready}) and expanding $(n \pm V/2)^{1+2/d} = n^{-1-2/d} (1 \pm V/2n)^{1+2/d}$ using the binomial theorem, all 
terms involving an odd power of $V$ cancel, and all terms involving an even power of $V$ are nonnegative, hence
	\[ \psi(A_n) - \psi(B_n) \geq C_2 V^2 n^{-1+2/d}. \]
Solving for $V$ and applying Proposition~\ref{mainprop} we obtain
	\begin{eqnarray*}  V &\leq& C_3 n^{1/2 - 1/d} (\psi(A_n) - \psi(B_n))^{1/2} \\
			     &\leq& C_3 n^{1/2-1/d} \Phi(n)^{1/2} + O(D^{1-1/8d} n^{1-1/2d}), \end{eqnarray*}
which yields (\ref{enough}). 
\end{proof}
 
\section{Isoperimetric Inequality for Exit Times}

In this section we show that among all regions $A \subset \Z^d$ of a given size, the maximal expected exit time
	\[ \bar{e}(A) = \sup_{x \in A} e_x(A) \] 
is asymptotically maximized when $A$ is a ball.

\begin{theorem}
\label{isoperimetric}
If $A \subset \Z^d$ and $|A|=n$, then
 	\begin{equation} \label{noname} \bar{e}(A) \leq e_o(B_n) + O(n^{2/d-\gamma_d} \log^2 n) \end{equation}
where 
	\begin{equation} \label{exponent} \gamma_d = \begin{cases} 1, & d=1 \\ \frac13, & d=2 \\ \frac{2^{-d}}{2d^2 \log 3}, & d\geq 3. 
\end{cases} \end{equation}
\end{theorem}

Recall from (\ref{martingale}) that $e_o(B_n) = \Theta(n^{2/d})$.  The case $d=1$ of Theorem~\ref{isoperimetric} is an elementary 
gambler's ruin calculation; for the remainder of this section we assume $d \geq 2$.  We have not attempted to optimize the exponent in 
the error term for $d\geq 3$, preferring instead to give the cleanest possible arguments.  A careful optimization of 
Lemma~\ref{orthant} should yield a somewhat smaller error.

Isoperimetric inequalities of this type have long been known for the exit time of Brownian motion from regions in $\R^d$.  The first such result goes back to P\'{o}lya 
\cite{Polya}, who showed that the disc in $\R^2$ maximizes ``torsional rigidity'' among all simply connected plane domains of a given area.  Aizenman and Simon 
\cite{AS} use a rearrangement inequality of Brascamp et al.\ \cite{BLL} to prove that a Euclidean ball in $\R^d$ simultaneously maximizes all moments of the Brownian exit 
time among all regions of a given volume.

The proof of Theorem~\ref{isoperimetric} will proceed in several steps.  We first appeal to a rearrangement inequality of Pruss 
\cite{Pruss} to reduce to the case when $A$ has a certain weak convexity property (Lemma~\ref{orthoconvex}).  This convexity enables us 
to estimate the exit time from points close to the boundary by bounding the hitting time of an orthant in $\Z^d$ (Lemma~\ref{orthant}).  
The Einmahl extension \cite{Einmahl} of the Koml\'{o}s-Major-Tusn\'{a}dy strong approximation \cite{KMT} (see 
also \cite{Zaitsev}) yields a close coupling of random walk in $A$ and Brownian motion in the corresponding region $A^\Box \subset 
\R^d$, so that the random walk is likely to be close to the boundary of $A$ when the Brownian motion exits $A^\Box$.  Finally, the 
theorem of Aizenman and Simon is used to bound the expected exit time of Brownian motion from $A^\Box$. 

Denote by $\varepsilon_1, \ldots, \varepsilon_d$ the standard basis for $\Z^d$, and by $H_i$ the hyperplane spanned by $\varepsilon_1, \ldots, \varepsilon_{i-1}, 
\varepsilon_{i+1}, \ldots, \varepsilon_n$.  Given a region $A \subset \Z^d$, for each $x \in H_i$ let
	\[ \alpha_i(x) = \# \{ j \in \Z ~|~ x+j\varepsilon_i \in A \}. \]
The {\it Steiner symmetrization} of $A$ with respect to the hyperplane $H_i$ is the region
	\[ \sigma_i A = \bigcup_{x \in H_i} \left\{ x+j\varepsilon_i \left| - \frac{\alpha_i(x)}{2} < j  \leq \frac{\alpha_i(x)}{2} 
\right. \right\}. \]
In words, $\sigma_i A$ is obtained by compressing to an interval each column of points in $A$ lying above a point $x \in H_i$, and then 
centering that interval about the hyperplane $H_i$, with preference for the positive side of the hyperplane if the interval has even 
length.  In particular, $|\sigma_i A| = |A|$.

We say that a region $A \subset \Z^d$ is {\it orthoconvex} if $x \in A$ and $x + k\varepsilon_i \in A$, $k > 0$ imply $x + j\varepsilon_i \in A$ for all $0 < j < k$; 
equivalently, any line in $\Z^d$ parallel to one of the coordinate axes meets $A$ in an interval (possibly empty).

\begin{lemma}
\label{orthoconvex}
For each $n \geq 1$ there exists an orthoconvex region $A \subset \Z^d$ which maximizes the quantity $\bar{e}(A)$ among all regions in 
$\Z^d$ of size $n$.
\end{lemma}

\begin{proof}
Denote by $\mathcal{A}$ the set of all connected regions $A \subset \Z^d$ of size $n$ containing the origin.  Clearly, the maximum 
value of $\bar{e}(A)$ among all regions of volume $n$ is attained by a region $A \in \mathcal{A}$.  It follows from the difference 
equation (\ref{laplacian}) and a rearrangement inequality of Pruss \cite{Pruss} that $\bar{e}$ does not decrease under Steiner 
symmetrization.  On the other hand, the quantity
	\[ \xi(A) := \sum_{x \in A} \sum_{i=1}^d |x_i+1/4| \]
strictly decreases under Steiner symmetrization unless $A$ is already symmetric.  Choosing from among those regions in $\mathcal{A}$ which maximize $\bar{e}$ one which 
minimizes $\xi$, we obtain a region that is Steiner symmetric about every coordinate axis, hence orthoconvex.
\end{proof} 

If $A$ is orthoconvex, then any point $x \in \partial A$ has a ``supporting orthant'' $Q$ based at $x$ lying entirely outside $A$.  To bound the time to exit $A$ from 
points near the boundary, it suffices to bound the time to hit this orthant.  We write ${\mathcal C}(x,r)$ for the $L^\infty$ ball of 
radius $r$ (cube of side length $2r+1$) centered at $x$.  Simple gambler's ruin considerations imply that
	\begin{equation} \label{ruin} E_x T_{\partial {\mathcal C}(x,r)} = O(r^2). \end{equation}
 
\begin{lemma} 
\label{orthant}
Let $Q$ be the nonnegative orthant $\{x \in \Z^d ~|~ x_i \geq 0, ~i=1,\ldots,d\}$, and let
	\[ p(k,r) = \sup_{||x||_\infty \leq k} \PP_x(T_Q>T_{\partial {\mathcal C}(o,r)}). \]
Then if $r \geq 3k$,
	\[ p(k,r) \leq \left( 1 - \frac{2^{-d}}{2d} \right) p(3k,r). \]
\end{lemma}

\begin{proof}
	\begin{figure}
	\centering
	\label{orthantfig}
	\includegraphics[scale=.75]{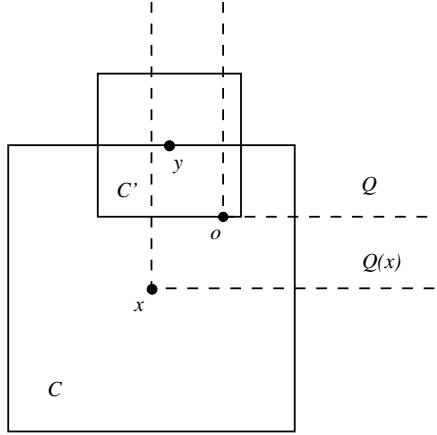}
	\caption{Diagram for the proof of Lemma~\ref{orthant}}
	\end{figure}
Given $x \in \Z^d$ with $||x||_\infty \leq k$, let $Q(x) = \{y \in \Z^d ~|~ y_i \geq x_i, ~i=1,\ldots,d\}$ be the orthant parallel to $Q$ based at $x$.  Subdividing the 
cube ${\mathcal C}={\mathcal C}(x,2k-1)$ into $2^d$ cubes of side length $2k-1$, the intersection $Q(x)\cap C$ consists of one of the 
cubes in the subdivision. By symmetry,
	\[ \PP_x(X_{T_{\partial C}} \in Q(x)) \geq 2^{-d}. \]
Now if $y$ is any point in $\partial C \cap Q(x)$, then an entire boundary face of the cube ${\mathcal C}' = {\mathcal C}(y,k-1)$ lies 
in $Q$ (Figure~2), hence by symmetry
	\[ \PP_y(X_{T_{\partial C'}} \in Q) \geq \frac{1}{2d}. \]
The result now follows from the observation that the $L^\infty$ norm of any point on the boundary of $C$ or $C'$ is at most $3k$.
\end{proof}

For a domain $D \subset \R^d$ denote by $T_{\partial D}^{\text{BM}}$ the time when
Brownian motion exits $D$.  The theorem of Aizenman and Simon \cite{AS} implies that $\EE_x T_{\partial D}^{\text{BM}}$ is maximized
among domains $D$ of volume $n$ when $D$ is a ball and $x$ is its center.  Since $\text{vol}(A^\Box)=|A|=n$ we obtain $\EE_x
T_{\partial A^\Box}^{\text{BM}} \leq (n/\omega_d)^{2/d}$ for all $x \in A^\Box$.  (We adopt the notational shorthand $\partial A^\Box 
:= \partial (A^\Box)$.)  Chebyshev's inequality gives
	\[ \sup_{x \in A^\Box} \PP_x(T_{\partial A^\Box}^{\text{BM}} > 2(n/\omega_d)^{2/d}) \leq \frac12, \]
hence 
	\begin{equation} \label{chebyshev} \PP_x(T_{\partial A^\Box}^{\text{BM}}>2m(n/\omega_d)^{2/d}) \leq 2^{-m}. \end{equation}
The following is a refinement of Lemma~\ref{orthant} in dimension two.

\begin{lemma}
\label{planarcase}
In dimension $d=2$, there are constants $a$ and $c$ such that
	\[ p(k,r) \leq c \left( \frac{k + a \log r}{r - a \log r} \right)^{2/3}. \]
\end{lemma}

\begin{proof}
Applying the map $z \mapsto z^{2/3}$, the conformal invariance of harmonic measure for planar Brownian motion implies \cite{LeGall} 
that
	\begin{equation} \label{2/3} p_{\text{BM}}(k,r) := \sup_{||x||_{\infty} \leq k} \PP_x(T_{\partial {\mathcal 
C}(o,r)^\Box}^{\text{BM}} < T_{Q^\Box}^{\text{BM}}) \leq c' \left( \frac{k}{r} \right)^{2/3}. \end{equation}
By the strong approximation theorem \cite{Einmahl,KMT,Zaitsev} there exists a constant $a>0$ and a coupling of simple random walk in 
$\Z^d$ with Brownian motion in $\R^d$, so that, except for an event $E_1$ of probability at most $\frac1n$, the coupled paths are 
separated by a distance of at most $a \log r-2$ up to time $s=\ceil{r^2 \log r}$.  Let $E_2$ be the event that the Brownian motion has 
not exited ${\mathcal C}(o,r)^\Box$ by time $s$.  By (\ref{chebyshev}) we have $\PP_x(E_2) = O(1/r)$.  On the event $E_1^c \cap E_2^c$, 
if the random walk exits  ${\mathcal C}(o,r)$ before hitting $Q$, the Brownian motion must exit ${\mathcal C}(o,r-a\log r)^\Box$ before 
hitting the translated quadrant $Q^\Box + (a \log r, a \log r)$; hence
	\[ p(k,r) \leq p_{\text{BM}}(k+a\log r,r-a\log r) + \PP_x(E_1) + \PP_x(E_2). \]
The result now follows from (\ref{2/3}).
\end{proof}

\begin{proof}[Proof of Theorem~\ref{isoperimetric}] 
Denote by $E_3$ the event that the random walk and Brownian motion paths in the strong approximation coupling are separated by distance 
more than $b \log n - d$ before time $s = \ceil{n^{2/d} \log n}$.  Choosing $b$ sufficiently large we can take $\PP(E_3)<\frac1n$.  
Write $T = T_{\partial A^\Box}^{\text{BM}}$, and denote by $E_4$ the event that $T>s$.  By (\ref{chebyshev}) we have $\PP_x(E_4) = 
O(\frac1n)$ for all $x \in A^\Box$.

On the event $E_3^c \cap E_4^c$ the location $X_T$ of the random walk when the Brownian motion exits $A$ is distance at most $b \log n$ 
from $\partial A$.  Let $Q \subset A^c$ be the supporting orthant at a point $Y \in \partial A$ within distance $b \log n$ of 
$X_T$.  For $j\geq 1$ let $F_j$ be the event that after time $T$ the walk travels to $L^\infty$ distance $3^j b \log n$ 
away from $Y$ before hitting $Q$.  Iteratively applying Lemma~\ref{orthant} with initial value $k = b \log n$ we obtain
	\begin{equation} \label{F_j} \PP_x(F_j) \leq \left( 1-\frac{2^{-d}}{2d} \right)^j \leq \exp \left( - \frac{2^{-d} j}{2d} 
\right). \end{equation}
Write $m = \ceil{\frac{\log n}{d \log 3}}$.  By (\ref{F_j}) we have
	\begin{equation} \label{F_m} \PP_x(F_m) = O(n^{-\gamma_d}), \qquad d \geq 3.  \end{equation}
In dimension two, Lemma~\ref{planarcase} with $k=b\log n$ and $r=3^j b \log n$ gives for $j \leq m$ 
	\begin{equation} \label{planarF_j} \PP_x(F_j) \leq C_0 3^{-2j/3}. \end{equation}
Taking $j=m$ we obtain $\PP_x(F_m) = O(n^{-1/3})$.  Thus (\ref{F_m}) holds in dimension two as well.
  

On the event $E_3^c \cap E_4^c \cap F_j^c$, the time for the random walk to exit $A$ is at most the exit time for Brownian motion plus
the time for the walk to go an additional $L^\infty$ distance $(3^j+1)b\log n$.  Hence
	\[ T_{\partial A} \leq T + \sum_{j=1}^m {\mathbb 1}_{F_j} T_{\partial 
{\mathcal C}(X_T,(3^j+1)b\log n)} + {\mathbb 1}_{F_m} \tilde{T} + {\mathbb 1}_{E_3}(s+\tilde{T_3}) + {\mathbb 
1}_{E_4}(s+\tilde{T_4}), \]
where $\tilde{T}$ is the additional time taken to exit $A$ if the walk travels to distance $3^m b \log n$ from $Y$ before hitting $Q$; 
and $\tilde{T_i}$ for $i=3,4$ is the additional time taken to exit $A$ after time $s$ if the event $E_i$ occurs.  Taking expectations 
and applying (\ref{ruin}), we obtain by the strong Markov property
	\begin{eqnarray} \label{expec} \EE_x T_{\partial A} &\leq& \left( \frac{n}{\omega_d} \right)^{2/d} + C_1 \sum_{j=1}^m 
\PP_x(F_j) 3^{2j} \log^2 n + O(n^{-\gamma_d}) \EE_x \tilde{T} + \frac{2s}{n} + \\
	~~~~~~~~~~&~&~~~~~~~~~~~~~~  + O \left(\frac1n \right) (\EE_x \tilde{T_3} + \EE_x \tilde{T_4}). \nonumber \end{eqnarray} 
By (\ref{F_j}), (\ref{F_m}) and (\ref{planarF_j}), 
	\[ \sum_{j=1}^m \PP_x(F_j) 3^{2j} \leq C_2 \PP_x(F_m) 3^{2m} = O(n^{2/d - \gamma_d} \log^2 n). \] 
Maximizing (\ref{expec}) over all $x \in A$ we obtain
	\[ \bar{e}(A) \leq \left( \frac{n}{\omega_d} \right)^{2/d} + O(n^{2/d-\gamma_d} \log^2 n) + O(n^{-\gamma_d}) \bar{e}(A) + 
O(n^{2/d-1} \log n) + \frac2n \bar{e}(A) \]
and solving for $\bar{e}(A)$ yields
	\[ \bar{e}(A) \leq \left(\frac{n}{\omega_d}\right)^{2/d} + O(n^{2/d - \gamma_d} \log^2 n). \qed \]
\renewcommand{\qedsymbol}{}
\end{proof} 

\section{Acknowledgements}
The authors thank Jim Propp for introducing the model, and for several helpful conversations.

\end{document}